\documentclass[a4paper]{article}
\usepackage[english]{babel}
\usepackage[utf8]{inputenc}
\usepackage{amsmath,amsthm}
\usepackage{amsfonts,amssymb}

\usepackage{tikz,tikz-cd} 

\usepackage{xifthen} 

\newcommand{\topos}{\textsc{Top}}

\newcommand{\Ob}{\mathbb{O}}

\newcommand{\Ecal}{\mathcal{E}} 
 
\newcommand{\Tcal}{\mathcal{T}}

\newcommand{\Ical}{\mathcal{I}} 
 
\newcommand{\Pcal}{\mathcal{P}} 
 
\newcommand{\Scal}{\mathcal{S}} 
 
\newcommand{\Fcal}{\mathcal{F}}

\newcommand{\Lcal}{\mathcal{L}}

\usepackage{titlesec}
\titleformat{\subsubsection}[runin]{\normalfont}{\thesubsubsection}{0pt}{}[.]

\renewcommand{\thesubsubsection}{\arabic{section}.\arabic{subsubsection}}
\csname @addtoreset\endcsname{subsubsection}{section}

\newcommand{\block}[1]
{

\par \subsubsection{} #1

\bigskip}

\newcommand{\blockn}[1]{\par #1 \bigskip}

\newcommand{\Prop}[1]
	{

	\bigskip
	
	\textbf{Proposition : }{\itshape #1}
		
	\bigskip
	
	}

\newcommand{\Cor}[1]
	{

	\bigskip
	
	\textbf{Corollary : }{\itshape #1}	
		
	\bigskip

	}

\newcommand{\Lem}[1]
	{

	\bigskip
	
	\textbf{Lemma : }{\itshape #1}
		
	\bigskip
	
	}

\newcommand{\Def}[1]
	{
	
	\bigskip
	
	\textbf{Definition : }{\itshape #1}
	
	\bigskip
	
	}

\newcommand{\Dem}[1]{
	
	\smallskip
	
	\textbf{Proof : } \par
	 {#1} $\square$
	 
	 \bigskip
}

\begin{document}

\pagestyle{plain}
\title{The localic isotropy group of a topos}
\author{Simon Henry}

\maketitle

\begin{abstract}

It has been shown by J.Funk, P.Hofstra and B.Steinberg that any Grothendieck topos $\Tcal$ is endowed with a canonical group object, called its isotropy group, which acts functorially on every object of the topos. We show that this group is in fact the group of points of a localic group object, called the localic isotropy group, which also acts on every objects, and in fact also on every internal locales and on every $\Tcal$-topos. This new localic isotropy group has better functoriality and stability property than the original version and shed some lights on the phenomenon of higher isotropy observed for the ordinary isotropy group. We prove in particular using a localic version of the isotropy quotient that any geometric morphism can be factored uniquely as a connected atomic geometric morphism followed by a so called ``essentially anisotropic'' geometric morphism, and that connected atomic morphism are exactly the quotient by open isotropy action.
\end{abstract}

\renewcommand{\thefootnote}{\fnsymbol{footnote}} 
\footnotetext{\emph{email:} simon.henry@college-de-france.fr}
\renewcommand{\thefootnote}{\arabic{footnote}} 


\tableofcontents

\bigskip

\section{Introduction}

\blockn{In \cite{funk2012isotropy}, J.Funk, P.Hofstra and B.Steinberg have introduced the idea of isotropy group of a topos. They have shown that any Grothendieck topos $\Tcal$ have a canonical group object $Z_{\Tcal}$ called the isotropy group of $\Tcal$ which acts (also canonically) on every object of $\Tcal$, and such that any morphism of $\Tcal$ is compatible with this action. They have also been considering the ``isotropy quotient'' $\Tcal_Z$ which is the full subcategory of $\Tcal$ of objects on which the action of $Z_{\Tcal}$ is trivial, it is a new Grothendieck topos (different from $\Tcal$ if $Z_{\Tcal}$ is non trivial) endowed with a connected atomic geometric morphism $\Tcal \rightarrow \Tcal_{Z}$. It also happen that in some case this topos $\Tcal_{Z}$ can have itself a non trivial isotropy group and this construction can be iterated, which has been referred to as ``higher isotropy'', but this does not happen in good case and seem to be somehow a pathological behavior.

\bigskip

It has also been observed that this isotropy group is the internal automorphism group of the universal point of $\Tcal$. For example if $\Tcal$ is a classifying topos $\Scal[\mathbb{T}]$ for some geometric theory $\mathbb{T}$, then the isotropy group is the internal automorphism group of the universal model of $\mathbb{T}$ in $\Scal[\mathbb{T}]$. This was first conjectured by Steve Awodey and a result of this kind appears in the PhD thesis of his student Spencer Breiner (\cite{breiner2014scheme}). From there, following some classical ideas from topos theory (see for example \cite{dubuc2003localic}) it is natural to look at the automorphism group of a point (or of a model of a theory) not as a discrete group but as a topological or better a localic group. This suggests that the isotropy group should arise naturally as a localic group.

\bigskip

The goal of this paper is to develop this idea: We introduce in section \ref{Sec_DefIsotropy} the localic isotropy group as the localic automorphism group of the universal point, which is described equivalently either by the fact that its classifies the data of a point of $\Tcal$ together with an automorphism of that point, or as the pullback the diagonal map of $\Tcal$ along itself. We show every object of the topos has a canonical action by the localic isotropy group and in fact more generally, every $\Tcal$-topos comes with such an action. Section \ref{Sec_DefIsotropy} also contains all the basic results about this isotropy group that does not really involve topos theory (mostly, the results of this section will be valid in any $2$-category with $(2,1)$-categorical finite limits). In section \ref{Sec_IsotropyQuotient} we introduce the notion of isotropy quotient adapted to the localic isotropy group, i.e. the fact that the subcategory of $\Tcal$ of object on which the isotropy action is trivial, is a topos $\Tcal_{G}$ endowed with a hyperconnected geometric morphism from $\Tcal \rightarrow \Tcal_{G}$. It is no longer the case in general that the quotient map is atomic. One can also consider isotropy quotient by arbitrary localic group endowed with an ``isotropy action'' i.e. a morphism to the isotropy group.

Section \ref{Sec_PostivieIsotropy} is the most important and technical. We focus on what happen when we take an isotropy quotient by a localic group which is locally positive (I.e. open, or overt), in this case one recover that the map $\Tcal \rightarrow \Tcal_G$ to the isotropy quotient is atomic and connected, and contrary to the ordinary case one get that the localic isotropy group of the isotropy quotient is nicely controlled by the isotropy group of the initial topos and the group which serve to construct the quotient, preventing in particular any higher isotropy phenomenon.

Conversely we also prove that any connected and atomic geometric morphism can be seen as an isotropy quotient by a locally positive isotropy group. Finally we see that any topos admit a ``maximal positive isotropy quotient'' which produces for any topos $\Tcal$ a connected atomic geometric morphism $\Tcal \rightarrow \Tcal_{I^+}$ where $\Tcal_{I^+}$ is ``essentially anisotropic'' i.e.  the isotropy group of $\Tcal_{I^+}$ have no locally positive sublocales other than ${1}$. Applying this to an arbitrary basis gives a unique factorization of any geometric morphism into a connected atomic morphism followed by an essentially anisotropic morphism, but this does not produces an orthogonal factorization system because the class of essentially anisotropic morphism is not stable under composition.

Finally in section \ref{Sec_Etale} we explain how the ordinary isotropy group mentioned in the beginning of the introduction (which we call the étale isotropy group) relate to our localic isotropy group and how the theory developed in \cite{funk2012isotropy} fits into ours.

}

\blockn{All the toposes are Grothendieck toposes over some base elementary topos $\Scal$ with a natural number object. By that we mean that they are (equivalent to) toposes of $\Scal$-valued sheaves over some $\Scal$-internal site, or equivalently that they are bounded $\Scal$-toposes.
Morphisms of toposes are the geometric morphisms over $\Scal$. The $2$-category of Grothendieck toposes and geometric morphisms over $\Scal$ is denoted $\topos$ (with the convention that $2$-morphisms are the natural transformation between the inverse image functors).

In particular everything done in this paper can be done over an arbitrary base topos and we will use this to obtain relative version of result proved over $\Scal$.
}

\section{The isotropy group}
\label{Sec_DefIsotropy}

\blockn{In this section we define the localic isotropy group of a topos $\Tcal$, its action on objects of $\Tcal$ as well as on $\Tcal$-topos and its basic properties. Despite being only formulated in terms of the category $\topos$ of bounded $\Scal$-toposes and its slices, the results of this section does not really uses much of topos theory, and most of it would hold in any weak $2$-category with finite limits (for the appropriate $(2,1)$-categorical notion of finite limits).}

\block{\Def{Let $\Tcal$ be a topos.

Let $\Ical_{\Tcal}$ be the contravariant weak functor from $\topos$ to the $2$-category of categories defined by:

\[\Ical_{\Tcal}(\Ecal) =\{ (f, \theta) | f:\Ecal \rightarrow \Tcal , \theta \text{ an automorphisms of $f^*$ } \} \]

morphisms in the category $\Ical_{\Tcal}(\Ecal)$ are the natural transformations $f^* \rightarrow f'^*$ compatible to the isomorphisms $\theta$ and $\theta'$.

If $v : \Ecal \rightarrow \Ecal'$ is a morphism of topos the functoriality $\Ical_{\Tcal}(v) : \Ical_{\Tcal}(\Ecal') \rightarrow \Ical_{\Tcal}(\Ecal)$ is given by pre-composition by $v$.

}

}

\block{

\Prop{The functor $\Ical_\Tcal$ is represented (up to equivalence) by a topos $I_{\Tcal}$ which can be defined as the $(2,1)$-categorical pullback:

\[
\begin{tikzcd}[ampersand replacement=\&]
I_{\Tcal} \arrow{r} \arrow{d} \& \Tcal \arrow{d}{\Delta} \\
\Tcal \arrow{r}{\Delta} \& \Tcal \times \Tcal \\
\end{tikzcd}
\]

Moreover the two arrows from $I_{\Tcal}$ to $\Tcal$ given by this diagram are isomorphic and corresponds to the map $(f,\theta) \mapsto f$ in terms of the representable functors.

}

\Dem{By the universal property of the pullback (and of the product) a morphism from $\Ecal$ to $I_{\Tcal}$ (defined as the above pullback) is given by two morphisms $(f,f')$ from $\Ecal$ to $\Tcal$ together with two isomorphisms $(\theta,\theta')$ between $f$ and $f'$.

The functors which send $(f,\theta)$ to $(f,f,\theta, Id_f)$ in one direction and $(f,g,\theta_1,\theta_2)$ to $(f, \theta_2^{-1}\circ \theta_1)$ in the other direction induces an equivalence of categories functorial in $\Ecal$ between the category of morphisms from $\Ecal$ to $I_{\Tcal}$ and the category $\Ical_{\Tcal}(\Ecal)$. 
}
}

\block{\Prop{ For any topos $\Tcal$, there is an internal localic group in $\Tcal$, also denoted $I_{\Tcal}$, such that the topos $I_{\Tcal}$ is the topos of $\Tcal$-valued sheaf over this internal localic group $I_{\Tcal}$. }

\Dem{By lemma 1.2 of \cite{johnstone1981factorization}, the fact that there is an internal locale $I_{\Tcal}$ in $\Tcal$ with this property corresponds to the fact that the geometric morphism $\Ical_{\Tcal} \rightarrow \Tcal$ constructed above is localic, which is the case because the map $\Delta: \Tcal \rightarrow \Tcal \times \Tcal$ is localic for any topos by \cite[B3.3.8]{sketches}, and that the pullback of a localic morphism is again localic by \cite[B3.3.6]{sketches}.

The fact that this internal locale $I_{\Tcal}$ has a group structure comes from the fact that the topos $\Ical_{\Tcal}$ has an obvious group structure over $\Tcal$ corresponding, in terms of the functor it represents to the composition of automorphisms of functors: the category of morphisms from a topos $\Ecal$ to $\Ical_{\Tcal} \times_{\Tcal} \Ical_{\Tcal}$ is equivalent to the category of morphisms from $\Ecal$ to $\Tcal$ endowed with two natural automorphisms, which can be composed or inverted functorially, providing this group structure.
}
}

\block{\Def{We call this internal localic group $I_{\Tcal}$ the localic isotropy group of $\Tcal$.}}

\block{Similarly to the isotropy group of \cite{funk2012isotropy}, the localic isotropy group $I_{\Tcal}$ is going to act on every object of $\Tcal$ making any morphism of $\Tcal$ $I_{\Tcal}$-equivariant, in fact it will acts on every $\Tcal$-topos in a way making the morphisms of $\Tcal$-topos $I_{\Tcal}$-equivariant. Before explaining this action, we would like to clarify a point about $2$-categories which will be central in the construction of this action.

Let $f: \Tcal \rightarrow \Ecal$ be a morphism of toposes, and assume that $\theta:f \rightarrow f$ is an automorphism of $f$. Then the pullback along $f$ functor has a natural automorphism induced by $\theta$ that we will denote $\Theta$. It can be represented as such: if $p:X \rightarrow \Ecal$ is a morphism, then the pullback along $f$ can be represented (in term of generalized points) as:

 \[ \{x \in X,t \in \Tcal, \gamma : p(x) \overset{\sim}{\rightarrow} f(t) \}, \]
 
and $\theta$ can also be represented in terms of generalized points as a functor which associate to each point $t \in \Tcal$ an automorphism $\theta_t : f(t) \overset{\sim}{\rightarrow} f(t)$, the automorphism of the pullback induced by $\theta$ can the be represented as:

\[\Theta_X : (x,t,\gamma) \mapsto (x,t,\theta_t \circ \gamma) \]
}

\blockn{
We are specifically interested in the case where $\Tcal = \Ecal$, $f$ is the identity of $\Tcal$, and $\theta$ is an automorphism of the identity of $\Tcal$. In this case the pullback along $f$ is (equivalent) to just the identity on the category of objects over $\Tcal$ so we should have an automorphism of every object ``over $\Tcal$'', but this is where there is a small $2$-categorical difficulty. Let us go though the description given above:

In terms of generalized point, $\Theta_X$ is the automorphism of $\{ x \in X, t \in \Tcal, \gamma: p(x) \overset{\sim}{\rightarrow} t \}$ which send $(x,t,\gamma)$ to $(x,t,\theta_x \circ \gamma)$. Hence the action is not on $X$ but on this object $X \times_{\Tcal} \Tcal$ which is isomorphic to $X$, the problem is that if one want to transport this action to an action on $X$ one get a trivial action: $\Theta_X$ is isomorphic to the identity of $X \times_{\Tcal} \Tcal$, indeed $Id : x \rightarrow x$  and $\theta_t : t \rightarrow t$ induces such an isomorphism from the identity to $\Theta_X$. But the point is that this is not a natural transformation over $\Tcal$ (it acts non trivially on the $t$ component) and hence this action can be non-trivial in the category $\topos_{/\Tcal}$. So if one want to express this action properly in terms of $X$, one needs $X$ to be ``fibrant''\footnote{This terms would only be completely accurate if ones where working in a strict $2$-category.} over $\Tcal$ (in a $2$-categorical sense) and this action corresponds to transport along $\theta_t$. In topos theoretical sense this is the case when $X$ is described by an internal site in $\Tcal$. 
}

\block{We now come back to our isotropy group. The morphism $p$ from $I_{\Tcal}$ to $\Tcal$, has an automorphism $\theta$ corresponding to the identity of $I_{\Tcal}$, it is the universal such automorphism in the sense that (by definition of $I_{\Tcal}$) for any automorphism of a morphism $\lambda$ of $\Ecal \rightarrow \Tcal$ there is a unique morphism $v$ from $\Ecal$ to $I_{\Tcal}$ such that $\lambda$ is obtained $\theta$ by composing with $v$.

This automorphism, along with the discussion above induces the action of $I_{\Tcal}$ on every $\Tcal$-topos: If $f:X \rightarrow \Tcal $ is a $\Tcal$-topos one can define an action $\alpha_X$ of $I_{\Tcal}$ on $X$ over $\Tcal$ as:

\[ \begin{array}{c c c}
 \{ i \in \Ical_{\Tcal} , x \in X, \gamma:f(x) \overset{\sim}{\rightarrow} p(i) \} & \rightarrow & \{t \in \Tcal, x \in X , \gamma:f(x) \overset{\sim}{\rightarrow} t \} \\
 (i,x,\gamma) & \mapsto & (p(i), x, i\circ \gamma) 
\end{array} 
\]

The first object is isomorphic to $I_{\Tcal} \times_{\Tcal} X$ and the second to $X$, but as above this map is trivial if one consider it as a map $I_{\Tcal} \times_{\Tcal} X \rightarrow X$, it is non trivial when considered as such a map over $\Tcal$.

\Prop{ $\alpha_X$ is an action of $I_{\Tcal}$ on $X$ for any $\Tcal$-topos $X$, any morphism $f: X \rightarrow X'$ over $\Tcal$ is equivariant for these actions.}

\Dem{This follows easily from the description given above.}
}

\block{Note that in the special case where $X$ is just an object of $\Tcal$ this action is a little simpler to describe: the pullback of $X$ to $I_{\Tcal}$ is endowed with an automorphism $\alpha_X$ induced by the canonical automorphism of the map $I_{\Tcal} \rightarrow \Tcal$, this automorphism is given by a map $I_{\Tcal} \times_{\Tcal} X \rightarrow X $ which can be checked to be an action of $I_{\Tcal}$. In the rest of the paper, we will mostly use the action of $I_{\Tcal}$ on objects, and we will anyway not need more than the action of $I_{\Tcal}$ on localic $\Tcal$-toposes, as localic $\Tcal$-toposes form an ordinary category, we do not really need to go in more details in these $2$-categorical problems.}

\block{One can also make a relative version of all of this:

\Def{If $f:\Ecal \rightarrow \Tcal$ is a $\Tcal$-topos, one defines the relative isotropy group $I_{\Ecal/\Tcal}$ as the isotropy group of $\Ecal$ when seen as a topos in $\Tcal$ instead of a topos over $\Scal$.

Equivalently,

\begin{itemize}

\item $I_{\Ecal/\Tcal}$ is the pullback:

\[
\begin{tikzcd}[ampersand replacement=\&]
I_{\Ecal/\Tcal} \arrow{r} \arrow{d} \& \Ecal \arrow{d}{\Delta} \\
\Ecal \arrow{r}{\Delta} \& \Ecal \times_{\Tcal} \Ecal \\
\end{tikzcd}
\]

\item A morphism to $I_{\Ecal/\Tcal}$ is the data of a morphism $v$ to $\Ecal$ together with an automorphism $\theta$ of $v$ such that the automorphism $f \circ \theta$ of $f \circ v$ is the identity.

\end{itemize}

}

One has a group homomorphism: $I_{\Ecal/\Tcal} \rightarrow I_{\Ecal}$ over $\Ecal$ given on generalized points by the obvious forgetful functor.

}

\block{If $G$ is any localic group over $\Tcal$ endowed with a morphism to $I_{\Tcal}$ then $G$ also acts on every object of $\Tcal$, as well as on every topos over $\Tcal$. Such a morphism from $G$ to $I_{\Tcal}$ will be called an \emph{isotropy action} of $G$. For example in the case of the natural morphism $I_{\Ecal/\Tcal} \rightarrow I_{\Ecal}$ one obtains an action of $I_{\Ecal/\Tcal}$ over every object of $\Ecal$ or topos over $\Ecal$. This action corresponds to the natural isotropy action seen internally in $\Tcal$ and the action is trivial on every object of $\Tcal$ or topos over $\Tcal$ when they are pulled back to $\Ecal$.}

\block{\Lem{Let $f :\Ecal \rightarrow \Tcal$ be a geometric morphism, then one has a pullback square:

\[
\begin{tikzcd}[ampersand replacement=\&]
f^{\sharp} I_{\Tcal} \arrow{r} \arrow{d} \& \Ecal \arrow{d}{\Delta} \\
\Ecal \times_{\Tcal} \Ecal \arrow{r}{Id \times Id} \& \Ecal \times \Ecal \\
\end{tikzcd}
\]

}
\Dem{In terms of generalized points, a morphism to the pullback is the data of: a morphism $g$ into $\Ecal$, two morphisms $h$ and $k$ into $\Ecal$, one isomorphism between $f \circ h$ and $f \circ k$, and two isomorphisms, one from $g$ to $h$, and one from $g$ and $k$. This is equivalent to one morphism $g$ into $\Ecal$ together with an automorphism of $f \circ g$ which is exactly a generalized point of $f^{\sharp} I_{\Tcal}$. This can also be alternatively explained in terms of pullback diagrams:

In
\[
\begin{tikzcd}[ampersand replacement=\&]
f^{\sharp} I_{\Tcal} \arrow{r} \arrow{d} \& I_{\Tcal} \arrow{d} \arrow{r} \& \Tcal \arrow{d}{\Delta} \\
\Ecal \arrow{r}{f} \& \Tcal \arrow{r}{\Delta} \& \Tcal \times \Tcal \\
\end{tikzcd}
\]

The two squares are pullback square hence the rectangle also is.

and in:

\[
\begin{tikzcd}[ampersand replacement=\&]
f^{\sharp} I_{\Tcal} \arrow{r} \arrow{d} \& \Ecal \times_{\Tcal} \Ecal \arrow{d} \arrow{r} \& \Tcal \arrow{d}{\Delta} \\
\Ecal \arrow{r}{\Delta} \& \Ecal \times \Ecal \arrow{r}{f \times f} \& \Tcal \times \Tcal \\
\end{tikzcd}
\]

the outer rectangle is a pullback because of the above remark and the rightmost square is a pullback by general property of fiber product, hence the leftmost square is a pullback which proves the lemma.

}

}

\block{\label{comparisonMap}\Prop{Let $f : \Ecal \rightarrow \Tcal$. There is a natural comparison map $I_{\Ecal} \rightarrow f^{\sharp} I_{\Tcal}$ which is a group morphism over $\Ecal$. Moreover this comparison maps fits into a pullback square:

\[
\begin{tikzcd}[ampersand replacement=\&]
I_{\Ecal} \arrow{r} \arrow{d} \& f^{\sharp} I_{\Tcal} \arrow{d} \\
\Ecal \arrow{r}{\Delta} \& \Ecal \times_{\Tcal} \Ecal \\
\end{tikzcd}
\]

}

\Dem{

The comparison map is easily defined in terms of generalized points:
A morphism into $f^*I_{\Tcal}$ is the data of a morphism $v$ into $\Ecal$ together with an automorphism of $f \circ v$.

A morphism into $I_{\Ecal}$ is the data of a morphism into $\Ecal$ together with an automorphism of this morphism. One can easily attached to it a morphism to $f^* I_{\Tcal}$ by simply applying $f \circ \_$ to the automorphism.

In terms of the pullback diagram of the lemma above, this comparison map appears to be the pullback:

\[
\begin{tikzcd}[ampersand replacement=\&]
I_{\Ecal} \arrow{r} \arrow{d} \& f^{\sharp} I_{\Tcal} \arrow{r} \arrow{d} \& \Ecal \arrow{d}{\Delta} \\
\Ecal \arrow{r}{\Delta} \& \Ecal \times_{\Tcal} \Ecal \arrow{r}{Id \times Id} \& \Ecal \times \Ecal \\
\end{tikzcd}
\]

Where the rightmost square is the pullback square of the lemma and the outer rectangle is the definition of $I_{\Ecal}$ as a pullback, which proves the existence of this comparison map, and moreover that this comparison map is a pullback of the diagonal map $\Ecal \rightarrow \Ecal \times_{\Tcal} \Ecal$. 
}
}

\block{\label{isotropyLeftexactsequence}\Prop{Let $f:\Ecal \rightarrow \Tcal$ be a geometric morphism, then $I_{\Ecal/\Tcal}$ is the kernel of the comparison map above, i.e. the sequence:

\[ 1 \rightarrow I_{\Ecal/\Tcal} \rightarrow I_{\Ecal} \rightarrow f^{\sharp} I_{\Tcal} \]

is exact.

In particular, if $I_{\Tcal} = \{1\}$ then $I_{\Ecal/\Tcal} \simeq I_{\Ecal}$.

}

\Dem{This is clear on generalized points: A morphism to $I_{\Ecal/\Tcal}$ is by definition a morphism $v$ to $\Ecal$ and an automorphism $\phi$ of $v$ such that $f \circ \phi$ is the identity. Such a couple $(v, \phi)$ without the last condition is the same as a morphism to $I_{\Ecal}$ and the last condition exactly says that the image into $f^{\sharp} I_{\Tcal}$ is the constant equal to the unit element.}

}

\section{Isotropy quotient}
\label{Sec_IsotropyQuotient}

\block{
\Def{If $G$ is any localic group over $\Tcal$ endowed with an isotropy action, i.e. a morphism to $I_{\Tcal}$, we define $\Tcal_{G}$ to be the full subcategory of $\Tcal$ of objects on which the isotropy action of $G$ is trivial. $\Tcal_G$ is called the isotropy quotient of $\Tcal$ by $G$.}

\Prop{$\Tcal_G$ is a topos, and the inclusion of $\Tcal_G$ in $\Tcal$ is the inverse image functor of a hyperconnected geometric morphism $p:\Tcal \rightarrow \Tcal_G$.}

See \cite[A4.6]{sketches} for the definition and basic properties of hyperconnected geometric morphisms.

\Dem{$\Tcal_G$ is a full subcategory of $\Tcal$ by definition, and because the action of $G$ is equivariant on all morphisms, it stable under ($\Scal$-indexed) colimits, finite limits, sub-objects and quotients. This is enough to imply the proposition.
}
}

\block{\label{Prop_FactoIsotrop}\Prop{Let $f:\Ecal \rightarrow \Tcal$ be a geometric morphism, let $G$ be a localic group over $\Ecal$ endowed with a morphism $v : G \rightarrow I_{\Ecal}$. The following conditions are equivalent:

\begin{itemize}

\item The composite $G \rightarrow I_{\Ecal} \rightarrow f^{\sharp}I_{\Tcal}$ is equal to $1$.

\item The morphism $v$ admit a (unique) factorization $G \rightarrow I_{\Ecal/\Tcal} \rightarrow I_{\Ecal}$.

\item The isotropy action of $G$ is trivial on every object of the form $f^*(X)$ for $X$ an object of $\Tcal$.

\item The geometric morphism $f$ factor into $\Ecal \rightarrow \Ecal_G \rightarrow \Tcal$

\end{itemize}

 }

\Dem{The equivalence of the first two points is exactly proposition \ref{isotropyLeftexactsequence}. The equivalence of the second and the third points follows immediately from the universal property of $I_{\Ecal/\Tcal}$, and the equivalence between the last two point follow immediately from the definition of $\Ecal_G$.} 
 
 }

\block{The proposition above has an important corollary:

\Cor{
Let $\Ecal \rightarrow \Tcal$ be an isotropy quotient of a topos $\Ecal$, I.e. $\Tcal = \Ecal_G$ for some localic group $G$ with an isotropy action $G \rightarrow I_{\Ecal}$, then $\Tcal = \Ecal_{I_{\Ecal/\Tcal}}$.

}

\Dem{ One has a factorization into $\Ecal \rightarrow \Ecal_{I_{\Ecal/\Tcal}} \rightarrow \Tcal$ corresponding to the relative isotropy quotient of $\Ecal$ over $\Tcal$.

And as $\Ecal_G$ factor into $\Tcal$ one has that $G \rightarrow \Ical_{\Ecal}$ factor into $I_{\Ecal/\Tcal}$ by proposition \ref{Prop_FactoIsotrop}, and hence a second factorization $\Ecal \rightarrow \Tcal \rightarrow \Ecal_{I_{\Ecal/\Tcal}}$, corresponding to the fact that $\Tcal$ is an isotropy quotient by a ``smaller'' isotropy action.

In both case the inverse image functor are inclusion of full subcategory so the existence of these two factorization implies the result.}

}

\block{The corollary above implies that one has a ``Galois theory'' classifying the isotropy quotient of a given topos $\Ecal$ in terms of certain subgroups of its isotropy group: those that arise has $I_{\Ecal/\Tcal}$ for some isotropy quotient $f:\Ecal \rightarrow \Tcal$. It is also not hard to see that any subgroup that appears as $I_{\Ecal/\Tcal}$ for $f :\Ecal \rightarrow \Tcal$ a general geometric morphism also appears as $I_{\Ecal/\Tcal'}$ for $\Tcal'$ the isotropy quotient $\Ecal_{I_{\Ecal\Tcal}}$. Unfortunately we are lacking of a good characterization of those.

\bigskip

\textbf{Open problem:} What are the subgroups of $I_{\Tcal}$ which appears as relative isotropy group $I_{\Tcal/\Ecal}$ of a geometric morphism $f: \Tcal \rightarrow \Ecal$ ?

}

\block{\label{WeirdExemple}Finally, it is important to note that without any assumptions on $G$ it is hard to say more about the map $\Tcal \rightarrow \Tcal_G$, here is an interesting example where this map is relatively general:

Let $\Tcal$ be the classifying topos of the theory of inhabited object, i.e. the theory with one sort $\Ob$, with no therms and with only one axioms: $\exists x \in \Ob$.

Equivalently, $\Tcal$ is the category of functors from the category of finite inhabited set to the category of sets, one takes $G$ be the full isotropy group of $\Tcal$, and we will see that the isotropy quotient $\Tcal_{G}$ is just the terminal topos, i.e. the category of sets.

\bigskip 

Indeed, for any finite inhabited set $X$, the functor $F \mapsto F(X)$ from $\Tcal$ to sets corresponds to a point of $\Tcal$, and its isomorphisms are exactly the isomorphisms of $X$, in particular, this shows that for $F$ in $\Tcal_G$ the action of the isomorphisms of $X$ on $F(X)$ should be trivial.

But one can easily that a presheaf satisfying this condition is automatically constant, and hence that $\Tcal_G$ is the category of sets.

Note that in this case, $\Tcal$ does not ``look like a category of group action'' at all and that the diagonal map $\Tcal \rightarrow \Tcal \times_{\Tcal_G} \Tcal$ is not a stable epimorphism (i.e. an epimorphism whose pullback are also epimorphisms) which is what we would need to apply the same kind of techniques as in the locally positive case treated in the next section.
}

\section{Locally positive isotropy}
\label{Sec_PostivieIsotropy}

\block{We start by some recall on local positivity and open maps:

\Def{An open subspace $U$ of a locale $X$ is said to be positive if whenever $U$ is written as a union of open subspaces:

\[ U = \bigcup_{i \in I} U_i \]

the indexing set is always inhabited: $ \exists i \in I$.

A locale is said to be locally positive if every open subspace can be covered by positive open subspaces.

}

If one uses classical logic, this notion is vacuous: ``positive'' is just equivalent to non-empty and every locale is locally positive, simply because any non-empty open subspace is the union of just itself and the empty open subspace is the union of the empty family. But within the internal logic a topos it is a non trivial notion:

\Prop{A locale $\Lcal$ internal to a topos $\Tcal$ is internally locally positive, if and only if the geometric morphism:

\[ sh_{\Tcal}(\Lcal) \rightarrow \Tcal \]

is open. It is an open surjection if and only if in addition $\Lcal$ is internally positive.
}

\Dem{This is \cite[C3.1.17]{sketches}}

So for a locale, locally positive is synonymous of ``open'' or ``overt''. We prefer the terminology ``locally positive'' to avoid the annoying double meaning of ``open sublocales''.
}

\blockn{Not that in a locally positive locale $\Lcal$ if $U = \bigcup U_i$ then one also have: \[ U =\bigcup_{i \text{ s.t. } U_i \text{ is positive.}} U_i . \]

Indeed, each $U_i$ is a union of positive open subspaces hence $U$ is the union of all the positive open subspaces which are included in one of the $U_i$, but each such open subspace is automatically included in a positive $U_i$, and hence  $U$ is the union of the positive $U_i$. }

\block{The main idea of this section, is that things works a lot better when one takes an isotropy quotient by a locally positive localic group rather than by a general group.}

\block{\label{QuotientByopengroup}\Lem{Let $G$ be a localic group over $\Tcal$ with a morphism $f : G \rightarrow I_{\Tcal}$, assume that the map from $G$ to $\Tcal$ is an open geometric morphism, then the map $p:\Tcal \rightarrow \Tcal_{G}$ is essential i.e. the inclusion functor $p^*$ has a left adjoint.}

\Dem{ Let $X$ be an object of $\Tcal$, and let $\theta_X :X \times G \rightarrow X$ be the action of $G$ on $X$. We will define an equivalence relation on $X$ by the following internal formula:

\[ x \sim y := \text{ the open subspace }\{ g \in G | gx =y \} \text{ is positive } \]

One can see that (Working internally in $\Tcal$ and using that internally $G$ is locally positive) it is an equivalence relation. Let $X_G$ be the quotient of $X$ by this relation. Then:

\begin{itemize}

\item The action of $G$ on $X_G$ is trivial, i.e. $X_G \in \Tcal_{G}$:

Indeed, as any map in $\Tcal$ the map quotient surjection $X \rightarrow X_G$ is $G$ equivariant. Internally, Let $x \in X$, $G$ is the union for $y \in X$ of the open subspace $G_{x,y} = \{ g | gx =y \}$. As $G$ is locally positive, one can also write $G$ as the union of those $G_{x,y}$ restricted to the $y$ such that $G_{x,y}$ is positive. In particular, all those $y$ are equivalent to $x$ and hence the action of $G$ on $x$ factor into the equivalence class of $x$ and hence is constant in $X_G$. This shows that $X_G$ is an object of $\Tcal_G$.

\item Every map from $X$ to an object of $\Tcal_{G}$ factor into $X_G$:

Let $f:X \rightarrow Y$ be any morphism, with $Y \in \Tcal_G$. Internally, Let $x,y \in X$. Internally in $G_{x,y}$, one can prove that $f(x)=f(y)$. Hence if $G_{x,y}$ is positive this implies that $f(x)=f(y)$ internally in $\Tcal$.

\end{itemize}

This shows that $X \mapsto X_G$ is a left adjoint to the forgetful functor $p^*:\Tcal_G \rightarrow \Tcal$ and hence concludes the proof.

}
}

\block{\label{AtomicQuotientMap}\Prop{If $G$ is a locally positive localic group in $\Tcal$ endowed with an isotropy action (i.e. in particular $f:G \rightarrow \Tcal$ is an open geometric morphism), then the natural map from $\Tcal \rightarrow \Tcal_G$ is atomic (and hyperconnected) in the sense of \cite[C3.5]{sketches}. }

\Dem{We will prove that the inclusion functor $f^*: sh(\Tcal_G) \rightarrow sh(\Tcal) $ is a logical functor, i.e. that it preserves the power object.

Let $X \in \Tcal_G$, let $\Pcal(X)$ be its power object in $\Tcal$, in order to see that $\Pcal(X)$ is also a power object in $\Tcal_G$ we just have to show that its natural $G$-action is trivial.

Let $Y$ be any object of $\Tcal$, and let $V \subset X \times Y$ be any sub-object.

$V$ is in particular stable under the action of $G$. In particular if, internally, $x,y \in Y$ are equivalent under the equivalence relation constructed in the proof of lemma \ref{QuotientByopengroup} and if $(v,x) \in V$ then $(v,y) \in V$ as well. Hence $V$ is the pullback of a subobject of $X \times Y_G$.

This proves that any morphism from $Y$ to $\Pcal(X)$ can be factored into a morphism from $Y_G$ to $\Pcal(X)$ and hence that $\Pcal(X)$ is in $\Tcal_G$ as claimed.
}
}

\block{Note that the conclusion of proposition \ref{AtomicQuotientMap} is obviously false without the local positivity assumption, in fact proposition \ref{AtomicMapAreGalois} below shows that an isotropy quotient $f:\Tcal \rightarrow \Ecal$ is locally positive only if it can be written as an isotropy quotient by a locally positive localic groups (although it may happen that a given isotropy quotient is both a quotient by a locally positive group and by a non locally positive group). The example given in \ref{WeirdExemple} also provides an explicit example where the map to the isotropy quotient is not atomic.}

\block{\label{atomicComparisonMap}\Prop{Let $f :\Ecal \rightarrow \Tcal$ be a connected atomic geometric morphism. Then the comparison map:

\[ I_{\Ecal} \rightarrow f^{\sharp} I_{\Tcal} \]

Is an open surjection.
}

So in some sense, in the case of a connected atomic morphism one has a short exact sequence of localic groups: $1 \rightarrow I_{\Ecal/\Tcal} \rightarrow I_{\Ecal} \rightarrow f^{\sharp} I_{\Tcal} \rightarrow 1$. Moreover as $f$ is an effective descent morphism of locales one can really think about $I_{\Tcal}$ as being the quotient of $I_{\Ecal}$ by $I_{\Ecal/\Tcal}$. This is the proposition that allows us to have some control on the isotropy group of the isotropy quotient in the case where the isotropy quotient is by a locally positive localic group. See proposition \ref{essentiallyAnisotropicQuotient} below for a typical examples of this sort of ideas.

\Dem{By proposition \ref{comparisonMap} the comparison map is a pullback of the diagonal map $\Ecal \rightarrow \Ecal \times_{\Tcal} \Ecal$. As open surjection are stable under pullback (see ), it is enough to show that the diagonal of a connected atomic topos is an open surjection. It is open because of \cite[C3.5.14]{sketches}, and it is a surjection by \cite[C3.5.6]{sketches} because it is a section of the morphism $\pi_1: \Ecal \times_{\Tcal} \Ecal \rightarrow \Ecal$ which is (hyper)connected and atomic by \cite[C3.5.12]{sketches}.
}
}

\blockn{We now want to show that among the locally positive localic group with an isotropy action there is a terminal object $I^+_{\Tcal}$ which defines a maximal connected atomic isotropy quotient. The idea is that thanks to the following proposition every locale (in particular $I_{\Tcal}$) has a maximal locally positive sublocales. In particular, the following proposition is specifically meant to be interpreted internally in a topos.}

\block{\Prop{Let $\Lcal$ be any locale, then:

\begin{itemize}

\item There is a maximal locally positive sublocale $\Lcal^+ \subset \Lcal$.

\item Any map from a locally positive locale to $\Lcal$ factor into the inclusion $\Lcal^+ \subset \Lcal$.

\item $\Lcal^+ \subset \Lcal$ is fiberwise closed (or weakly closed, see \cite{sketches} just before C1.1.22) inside $\Lcal$.

\item If $G$ is a localic group then $G^+$ is a localic subgroup

\end{itemize}

}

\Dem{The existence of $\Lcal^+$ follows from the fact that a co-product of a small family of locally positive locales is again locally positive and that if $X$ is locally positive and $f:X \rightarrow \Lcal$ is a morphism then the (regular) image of $X$ in $\Lcal$ is locally positive. This also implies the second point. The third point follows from the fact that the fiberwise closure of $\Lcal^+$ in $\Lcal$ is itself locally positive by \cite[C3.1.14(ii)]{sketches}. As for The last point: the terminal locale is locally positive, hence the unit of $G$ lies in $G^+$, and as $G^+$ and $G^+ \times G^+$ are both locally positive, the inversion and the multiplication map restrict as maps $G^+ \rightarrow G^+$ and $G^+ \times G^+ \rightarrow G^+$, and hence $G^+$ is a subgroup.
}

}

\block{\Def{One says that a geometric morphism $f : \Ecal \rightarrow \Tcal$ is completely anisotropic if $I_{\Ecal/\Tcal} = \{1\}$ and essentially anisotropic if $I_{\Ecal/\Tcal}^+ = \{1\}$.} 

}

\block{\label{essentiallyAnisotropicQuotient}\Prop{Let $f:\Ecal \rightarrow \Tcal$ be a geometric morphism, let $G$ be $(I_{\Ecal/\Tcal})^+$ endowed with its natural inclusion map to $I_{\Ecal/\Tcal}$, then the geometric morphism $\Ecal_G \rightarrow \Tcal$ is essentially anisotropic.}

\Dem{To simplify notation, all the isotropy groups are considered over $\Tcal$. Let $p$ be the map $\Ecal \rightarrow \Ecal_G$. It is connected and atomic by proposition \ref{atomicComparisonMap} because $G$ is locally positive.

We want to prove that $I_{\Ecal_G}^+=\{1\}$, i.e. that any morphism from a locally positive locale $X$ to $I_{\Ecal_G}$ is constant equal to the unit element. As the map $p: \Ecal \rightarrow \Ecal_G$ is hyperconnected, it is in particular a stable surjection, hence it is enough to prove that any map over $\Ecal$ from a locally positive $\Ecal$-locale $X$ to $p^{\sharp} I_{\Ecal_G}$ is constant equal to the unit element. We fix such a map.

As the comparison map from $I_{\Ecal}$ to $p^{\sharp} I_{\Ecal_G}$ is an open surjection (by \ref{atomicComparisonMap}), if one form the pullback $Y = X \times_{p^{\sharp} I_{\Ecal_G}} I_{\Ecal}$ then the projection $Y \rightarrow X$ is also an open surjection, hence $Y$ is locally positive, and hence the second projection $Y \rightarrow I_{\Ecal}$ factor into $I_{\Ecal}^+ = G$.

But $G$ is in the kernel of the comparison map $I_{\Ecal} \rightarrow p^{\sharp} I_{\Ecal_G}$ hence, as $Y \rightarrow X$ is an open surjection, this implies that the map from $X$ to $I_{\Ecal_G}$ is constant equal to $1$ and hence proves the result.
}
}

\blockn{We are now ready to prove that conversely any connected atomic map is canonically an isotropy quotient by a locally positive isotropy group:}

\block{\label{AtomicMapAreGalois}\Prop{Let $ f : \Ecal \rightarrow \Tcal$ be a connected atomic morphism, then:

\begin{itemize}

\item The relative isotropy group $I_{\Ecal/\Tcal}$ is locally positive in $\Ecal$.

\item The topos $\Ecal \times_{\Tcal} \Ecal$ is equivalent to the topos of object of $\Ecal$ endowed with a $I_{\Ecal/\Tcal}$-action. Under this identification, $\Delta^*$ is the functor that forget the action, $\pi_1^*$ is the functor that endows an object with the trivial action and $\pi_2^*$ is the functor that endows an object with its canonical $I_{\Tcal/\Ecal}$-action.

\item The natural map $\Ecal_{I_{\Ecal/\Tcal}} \rightarrow \Tcal$ is an equivalence of topos.

\end{itemize}

}

For the proof of this proposition we will need to use some results from descent theory. We refer the reader to \cite[C5.1]{sketches} for an introduction to descent theory which contains already a lot more than what we need.

\Dem{

One has a pullback square:

\[
\begin{tikzcd}[ampersand replacement=\&]
I_{\Ecal/\Tcal} \arrow{r} \arrow{d} \& \Ecal \arrow{d}{\Delta} \\
\Ecal \arrow{r}{\Delta} \& \Ecal \times_{\Tcal} \Ecal \\
\end{tikzcd}
\]

but the arrow $\Ecal \rightarrow \Ecal \times_{\Tcal} \Ecal$ is open, hence $I_{\Ecal/\Tcal} \rightarrow \Ecal$ also is, which proves the first point.

The map $\pi_1 : \Ecal \times_{\Tcal} \Ecal \rightarrow \Ecal$ corresponds internally in $\Ecal$ to a connected atomic topos which has a point given by $\Delta:\Ecal \rightarrow \Ecal \times_{\Tcal} \Ecal$ hence by \cite[C5.2.13]{sketches} it can be identified with the topos of objects of $\Ecal$ endowed with an action of the localic automorphism group of $\Delta$, but this is (by definition) the isotropy group $I_{\Ecal/\Tcal}$. Following the construction of the isotropy action shows that $\pi_2^*$ indeed corresponds to endowing objects with its isotropy action.

Moreover, $f$ (as any hyperconnected morphism) is an effective descent morphism for objects. Hence $\Tcal$ is equivalent to the category of objects of $\Ecal$ endowed with a descent data relative to $p$, once we replace $\Ecal \times_{\Tcal} \Ecal$ by the topos of objects of $\Ecal$ endowed with an action of $I_{\Ecal/\Tcal}$, this descent data is described as an isomorphism between an object $X$ with the trivial isotropy action and $X$ with the canonical isotropy action which is the identity on $X$, hence the category of such objects endowed with a descent data is just the category of objects whose isotropy action is trivial, which proves the third point.
}

The exact same proof, together with the fact that hyperconnected morphisms are also effective descent morphisms for locales, actually proves a stronger result: the category of locales over $\Tcal$ is equivalent to the full subcategory of locales over $\Ecal$ which have a trivial isotropy action.

In fact, using I.Moerdijk's result from \cite{moerdijk1989descent} that open surjections (in particular hyperconnected morphisms) are effective descent morphisms in the $2$-category of toposes, one can even deduce that the category $\topos_{/\Tcal}$ is the equivalent to the category of toposes over $\Ecal$ endowed with a trivialization over $\Ecal$ of their isotropy action.
}

\block{\Prop{Any geometric morphism $f: \Ecal \rightarrow \Tcal$ has a unique factorization (up to unique isomorphisms) as a connected atomic morphism followed by an essentially anisotropic morphism given by:

\[ \Ecal \rightarrow \Ecal_{I_{\Ecal/\Tcal}^+} \rightarrow \Tcal \]

}

\Dem{The factorization given in the proposition is clearly a factorization as an atomic connected morphism (by proposition \ref{essentiallyAnisotropicQuotient}) followed by an essentially anisotropic morphism (by proposition \ref{atomicComparisonMap}). We will now prove the uniqueness of the factorization.

Let $\Ecal \overset{p}{\rightarrow} \Fcal \overset{a}{\rightarrow} \Tcal$ be any such factorization. Let $G = I_{\Ecal/\Fcal}$, by proposition \ref{AtomicMapAreGalois}, $G$ is locally positive and $\Fcal$ is canonically isomorphic to $\Ecal_{G}$, and by propositions \ref{isotropyLeftexactsequence} and \ref{atomicComparisonMap} applied with $\Tcal$ as a basis one has:

\[ 1 \rightarrow G \rightarrow I_{\Ecal/\Tcal} \twoheadrightarrow p^{\sharp} I_{\Fcal/\Tcal} \rightarrow 1 \]

$G \rightarrow I_{\Ecal/\Tcal}$ factor into $I_{\Ecal/\Tcal}^+$ because $G$ is locally positive. Over $\Fcal$, as $p$ is an open map, the topos $I_{\Ecal/\Tcal}^+ \rightarrow \Ecal \rightarrow \Fcal$ is open, and hence its map to $I_{\Fcal/\Tcal}$ has to factor into $I_{\Fcal/\Tcal}^+$ which is $\{1\}$ because one has assumed that $\Fcal \rightarrow \Tcal$ is essentially anisotropic. Hence the map $I_{\Ecal/\Tcal}^+ \rightarrow p^{\sharp} I_{\Fcal/\Tcal}$ is constant, and hence $G=I_{\Ecal/\Tcal}^+$ which concludes the proof.
}

}

\block{\label{exemplEtendu}One does not get an orthogonal factorization system or the unique lifting property because essentially anisotropic map are not stable under composition:

Let $\Tcal$ be the topos of sheaves over $[-1,1]$ equivariant for the natural multiplication action of $\{1,-1\}$, and Let $\Ecal$ be the topos of set endowed with an action of $\{-1,1\}$. There is a geometric morphism from $\Ecal$ to $\Tcal$ whose inverse image functor is the germ at $0$ with the induced action of $\{-1,1\}$. The topos $\Tcal$ is essentially anisotropic (its non trivial isotropy is concentrated over the subspace $\{0 \} \subset [-1,1]$ of empty interior and hence cannot contains any locally positive locales), the morphism from $\Ecal$ to $\Tcal$ is localic hence completely anisotropic but the topos $\Ecal$ is not completely anisotropic.

On the other hand, completely anisotropic maps are stable under composition because of proposition \ref{isotropyLeftexactsequence} applied relatively to the target of the composition, but it is not clear at all that is produces a factorization system as general isotropy quotient by non locally positive group can be relatively wild and we do not know if for example the isotropy quotient by the full isotropy group is always completely anisotropic or not.
}

\section{Comparison to the ``étale isotropy group''}
\label{Sec_Etale}

\block{\Prop{The isotropy group $Z_{\Tcal}$ of $\Tcal$, as defined in \cite{funk2012isotropy}, is the group of points of $I_{\Tcal}$.}

We will call $Z_{\Tcal}$ the étale isotropy group.

\Dem{Let $Z_{\Tcal}$ be the group of points of $I_{\Tcal}$. For any object $X$ of $\Tcal$, the morphisms from $X$ to $Z_{\Tcal}$ are the same as the morphisms of toposes over $\Tcal$ from $\Tcal_{/X}$ to $I_{\Tcal}$, hence they are the same as automorphisms of the morphism from $\Tcal_{/X}$ to $\Tcal$, which is exactly the universal property of the isotropy group defined in \cite{funk2012isotropy} (the group structure and the functoriality are immediately checked to be the same).}
}

\blockn{Note that the étale isotropy group, as it is étale, is always locally positive, so the isotropy quotient constructed in \cite{funk2012isotropy} fits into the theory of the previous section.}

\block{This localic pictures ``explains'' the higher isotropy phenomenon of \cite{funk2012isotropy}:

We start with a topos $\Tcal$, and $I$ its localic isotropy group. Taking the isotropy quotient in the sense of \cite{funk2012isotropy} amount to taking the isotropy quotient by $Z \rightarrow I$. The resulting topos has an isotropy group which is (up to descent) a quotient $I_1$ of $I$ by at least the group of all the points of $I$. But it is not because all the points of $I$ have been killed in $I_1$ that $I_1$ cannot have some new points that does not lift into points of $I$. Hence the isotropy quotient can still have a non trivial étale isotropy group, and this phenomenon cannot happen when one quotient directly by the full locally positive part of the isotropy group as in the previous section.
}

\block{Finally, there is one case where the two theories agrees:

\Prop{If the unit map $\Tcal \rightarrow I_{\Tcal}^+$ is open, then $I_{\Tcal}^+$ is discrete and is the étale isotropy group $Z_{\Tcal}$. This happens if $\Tcal$ is locally essentially anisotropic, or for example if it is an étendu. }

\Dem{If the unit map of a localic group is open, then the diagonal map $G \rightarrow G \times G$ of the group is also open (because it is the pullback of the unit map along the multiplication map). The group of points always factor into $I_{\Tcal}^+$, but under the assumption of the proposition $I_{\Tcal}^+$ ends up being locally positive with an open diagonal, hence is discrete by \cite[C3.1.15]{sketches} and hence it is exactly the group of points.

For any slice $\Tcal_{/X}$ of  $\Tcal$  the comparison map:

\[ I_{\Tcal_{/X}} \rightarrow X \times_{\Tcal} I_{\Tcal} \]

is an open inclusion because it is a pullback of the diagonal map $\Tcal_{/X} \rightarrow \Tcal_{/X} \times_{\Tcal} \Tcal_{/X}$, which is an open inclusion.

Note that if $\Lcal$ is a locale over $\Tcal_{/X}$ then $\Lcal^+$ is the same whether one see $\Lcal$ as a locale in $\Tcal_{/X}$ or as a locale in $\Tcal$ with a map to $X$, and it corresponds (internally in $\Tcal$) to apply $^+$ to every fibers of this map to $X$.

In particular one has that
\[ (I_{\Tcal} \times X)^+ =  I_{\Tcal}^+ \times X, \]

and $I_{\Tcal_{/X}}^+ = (I_{\Tcal} \times X)^+ \cap I_{\Tcal_{/X}}$.

hence if there is a $X$ is such that  $I_{\Tcal_{/X}}^+ = \{1\}$, this proves (assuming $X$ is inhabited) that the map $1 \rightarrow I_{\Tcal}^+$ is open.

}

}

\block{Remarks: Even for étendu, the localic group $I_{\Tcal}$ does not have to be discrete. The example mentioned in \ref{exemplEtendu} of the topos of equivariant sheaves over $[-1,1]$ with the action of $\{-1,+1 \}$ by multiplication, is an étendu with non trivial isotropy group (because the point corresponding to $0$ has a non trivial automorphism) but the étale isotropy group and the positive isotropy group (isomorphic because of the proposition above) are trivial as mentioned earlier.}

\block{Finally, while the localic theory explains and somehow solve the higher isotropy phenomenon observed in \cite{funk2012isotropy}, it is not clear that it does not produce a new sort of ``higher isotropy''. More precisely, we have very little control on the isotropy quotient by an isotropy group which is not locally positive\footnote{One can also probably develop a similar theory to control isotropy quotient by \emph{compact} localic groups, with no local positivity assumption.}, and we do not know the answer to the following question:

\bigskip

\textbf{Open problem:} Given $\Tcal$ a topos and $G$ its full localic isotropy group can the isotropy quotient $\Tcal_{G}$ have a non trivial localic isotropy group ?

}

\bibliographystyle{plain}

\end{document}